\documentclass[12pt]{article}

\usepackage{t-angles}
\usepackage[all]{xy}
\usepackage{amsmath,amsthm, amscd, amssymb, amsfonts}

\usepackage{hhline}

\title{ Quantum Theory in Design }
\author{  \small  Zhengtang Tan,  Shouchuan Zhang \\
\small Department  of Mathematics, Hunan University\\
\small   Changsha  410082, \ P.R. China,
\small Emails: z9491@sina.cn  }
\date{}

\begin{document}
% \tableofcontents
\newtheorem{Theorem}{\quad Theorem}[section]
\newtheorem{Proposition}[Theorem]{\quad Proposition}
\newtheorem{Definition}[Theorem]{\quad Definition}
\newtheorem{Corollary}[Theorem]{\quad Corollary}
\newtheorem{Lemma}[Theorem]{\quad Lemma}
\newtheorem{Example}[Theorem]{\quad Example}
\newtheorem{Remark}[Theorem]{\quad Remark}

\numberwithin{equation}{section}

\maketitle \addtocounter{section}{-1}

\numberwithin{equation}{section}

\begin {abstract}

In this paper, we use braiding diagrams to present rules of shapes and designs. That is, we design colour, design size, design brightness, design codes by means of braiding.

\vskip.2in
2000 Mathematics Subject Classification: 16B50.

Keywords: code,  braiding,  design.
\end {abstract}

%%%%%%%%%%%%%%%%%%%%%%%%%%%%%%%%%%%%%%%%%%%%%%%%%%%%%%%%%%%%%%%%%%%%%%%%%%%%%%
\section{Introduction}\label {s0}
%%%%%%%%%%%%%%%%%%%%%%%%%%%%%%%%%%%%%%%%%%%%%%%%%%%%%%%%%%%%%%%%%%%%%%%%%%%%%%

The Yang-Baxter equation first came up in a paper by Yang as
factorization condition of the scattering S-matrix in the many-body
problem in one dimension and in the work of Baxter on exactly
solvable models in statistical mechanics. It has been playing an
important role in mathematics and physics ( see \cite {BD82} , \cite
{YG89} ). Attempts to find solutions of The Yang-Baxter equation in
a systematic way have led to the theory of quantum groups \cite {Ka95}.

Shape grammars have been used as a computational design tool  for over two decades. Shape grammars are a
production system created by taking a sample of the whole for which one
is trying to write a language \cite {St06}. From this sample a vocabulary of shapes
can be written that represent all the basic forms of that sample. By defining
the spatial relationships between those forms and how the forms are related
to each other, shape rules can be written. A shape rule consists of a
left and right side. If the shape in the left side matches a shape in a drawing
then the rule can be applied, and the matching shape changes to match
the right side of the rule. The shape rules allow the addition and  subtraction
of shapes, which in the end are perceived as shape modifications.
These shape rules, combined with an initial shape, produce a shape grammar
that represents the language of the design \cite {St06}. Shapes themselves can
exist as points, lines, planes, volumes, or any combination thereof \cite {St06}. All
shape generation must start with an initial shape: a point, a coordinate axis,
or some foundation from which to start the shape grammar. If the grammar
is going to end, it can end with a terminal rule, which prevents any
other rules from being applied after it. This forces there to be closure in
the rule sequence. Alternatively, a design sequence can continue indefinitely
and designs could be chosen at any point in the design process.
The method discussed here fundamentally changes the method of developing
the shape grammar. Mackenzie \cite {Ma69} demonstrates in a simplified
case that if the fundamental shapes in the language are defined, the relationship
between the shapes can be inferred through examples. These inferred
relationships can then be mapped to trees which in turn can be used
to automatically create shape grammar rules for that language. While this
was shown effective in a particular use, it is more common in practice that
the vocabulary of the design to be used as the foundation of a shape
grammar is determined by the creator of the grammar. The creator of the
grammar looks at the sample and subjectively derive the vocabulary.
From that vocabulary, the rules are formed based upon the creator¡¯s experience
and intention. It is quite possible that two different persons looking
at the same sample of shapes would create two very different shape grammars.

  In \cite   {OCB08}, the results from the principal component
analysis are used to create a new coupe shape grammar based upon these
discovered shape relationships. The shape grammar is then used to create
new coupe vehicles. Although the focus of this paper is on vehicle design,
the methods developed here are applicable to any class of physical products
based on a consistent form language. In \cite {Ho10}, visual 3D spatial grammars are studied.

In this paper we use braiding diagrams to represent rule of shapes and designs. This simplifies the  expression
 of designs. The shape rules allow  braiding of two  shapes, which in the end are perceived as shape modifications. We represent new design by means of  restricted PBW of Nichols algebras. That is, we   designs in  braided tensor categories.  We design colour, design size, design brightness, design codes by means of braiding.
 This completely
is a new method and will be applied more and more in various  designs.

\section {Preliminaries} \label {c1}

We begin with the tensor category (see \cite {Ma98}, \cite {Ka95} and \cite {Zh99}). We define the product ${\cal C}
\times {\cal D}$ of two category ${\cal C}$  and ${\cal D}$ whose
objects are pairs of objects $(U,V) \in (ob {\cal C}, ob {\cal D})$
and whose morphisms are given by
$$Hom _{{\cal C}\times {\cal D}}((V,W) (V',W'))  = Hom _{\cal C} (V,V')
\times Hom _{\cal D} (W, W'). $$ Let ${\cal C}$  be a category and
$\otimes $   be a  functor from ${\cal C} \times {\cal C} $  to
${\cal C}$. This means

(i) we have object $V \otimes W$  for any $V, W \in ob {\cal C};$

(ii) we have morphism $f \otimes g$  from $U\otimes V$  to $X\otimes
Y$ for any morphisms $f$  and $g$  from $U$ to $X$  and from $V$  to
$Y$;

(iii)  we have
$$(f \otimes g) (f' \otimes g') = (ff') \otimes (gg') $$
for any morphisms $f: U\rightarrow X, g : V\rightarrow Y, f' :
U'\rightarrow V$ and  $g' : V' \rightarrow V;$

(iv) $id_{U\otimes V} = id_U \otimes id _V.$

Let $\otimes \tau $ denote the functor from ${\cal C} \times {\cal
C} $  to ${\cal C}$ such that $(\otimes \tau)(U, V) =( V \otimes U)$
and  $(\otimes \tau )(f, g) = g \otimes f$,
 for any objects $U, V, X, Y$ in ${\cal C},$
and for any morphisms $f : U \rightarrow X$ and $g : V \rightarrow
Y.$

 An associativity constraint
$a$ for tensor $\otimes $ is a natural isomorphism
$$a: \otimes (\otimes \times id ) \rightarrow \otimes (id \times \otimes ).$$
This means that, for any triple $(U,V,W)$  of  objects of ${\cal
C}$, there is a morphism $a_{U,V,W} : (U\otimes V) \otimes W
\rightarrow U\otimes (V\otimes W)$  such that
$$a_{U',V',W'} ((f\otimes g) \otimes h) = (f\otimes (g\otimes h) ) a_{U,V,W} $$
for any morphisms $f, g$ and $h$  from $U$ to $U'$, from $V$  to
$V'$ and from $W$ to $W'$ respectively.

Let $I$ be an object of ${\cal C}$. If there exist  natural
isomorphisms
  $$ l : \otimes (I \times id )  \rightarrow  id  \hbox { \ \ \ and \ \ \ }
   r : \otimes (id \times I )  \rightarrow  id \ , $$
   then $I$  is called the unit object of ${\cal C}$
   with left unit constraint $l$
   and right unit constraint $r$.

\begin {Definition}  \label {0.1}
$({\cal C}, \otimes, I,  a, l, r )$ is called a tensor category if
${\cal C}$  is equipped with a tensor product $\otimes $, with a
unit object $I$, an associativity constraint   $a$ , a left unit
constraint $l$  and a right unit constraint $r$ such that the
Pentagon Axiom and the Triangle Axiom are satisfied, i.e.
$$(id _U \otimes a_{V,W,X}) a_{U, V\otimes W, X} (a_{U,V,W} \otimes
id _X) = a_{U,V,W\otimes X}a_{U\otimes V,W,X}$$ and
$$(id _V \otimes l_W) a_{V,I,W} = r_{V} \otimes id _W$$
for any $U, V, W, X \in ob {\cal C}.$

Furthermore, if  there exists a  natural isomorphism
$$C: \otimes \rightarrow \otimes \tau$$
such that the Hexagon Axiom holds, i.e.

$$a_{V,W,U} C_{U,V\otimes W} a_{U,V,W} =
(id _V \otimes C_{U,W}) a_{V,U,W}(C_{U,V}\otimes id _W)$$ and
$$a^{-1}_{W,U,V} C_{U\otimes V, W} a^{-1}_{U,V,W} =
( C_{U,W} \otimes id_ V ) a^{-1}_{U,W,V}(id_U \otimes C_{V,W}),$$
for any $U,V,W \in ob {\cal C},$ then  $({\cal C}, \otimes, I,  a,
l, r, C )$ is called a braided
 tensor category.
In this case, $C$ is called a braiding of ${\cal C}.$  If $C_{U,V} =
C_{V,U} ^{-1}$ for any $U, V \in ob {\cal C},$ then $({\cal C}, C)$
is called a symmetric braided tensor category or a symmetric tensor
category. Here functor $\tau: {\cal C} \times {\cal C} \rightarrow
{\cal C} \times {\cal C}$ is the flip functor defined by $\tau
(U\times V) = V\times U$ and $\tau (f \times g) = g \times f,$ for
any $U,V \in {\cal C},$  and morphisms $f$ and $g.$ Note that we
denote the braiding $C$ in braided tensor category $({\cal C},
\otimes, I,  a, l, r, C )$
 by $^{\cal C}C$ sometimes.

\end {Definition}

\begin {Example} \label {0.4} (The tensor category of vector spaces )
The most fundamental   example of a tensor category is given by the
category ${\cal C} = {\cal V}ect (k)$ of
 vector
spaces over  field $k$. ${\cal V}ect(k)$  is equipped with tensor
structure for which $\otimes $ is the tensor product of the vector
spaces over $k$, the unit object $I$ is the ground field $k$ itself,
and the associativity constraint and unit constraint are the natural
isomorphisms
$$a_{U, V, W}((u \otimes v) \otimes w):= u \otimes (v \otimes w)  \hbox { \ \ and \ \ }
l_V(1 \otimes v) := v := r_V (v \otimes 1)$$ for any vector space $U,
V, W$ and $u \in U, v\in V, w\in W.$

Furthermore, the most fundamental example of a braided tensor
category is given by the tensor category  ${\cal V}ect (k)$, whose
braiding is usual twist map from $U \otimes V $  to $V\otimes U$
defined by sending $a \otimes b$ to $b \otimes a$ for any $a\in U,
b\in V$.
\end {Example}

Now we define some notations. If $f$ is a morphism from $U$ to $V$
and $g$  is a morphism from $V$ to $W$, we denote the composition
$gf$ by
\[
\begin{tangle}
\object{U}\\
\vstr{50}\id\\
\O {gf}\\
\vstr{50}\id\\
\object{W}\\
\end{tangle}
\step=\step
\begin{tangle}
\object{U}\\
\O {f}\\
\O {g}\\
\object{W}\\
\end{tangle}\ \ \ .
\]

We usually  omit $I$ and  the morphism $id _I$ in any diagrams. In
particular, If $f$ is a morphism from $I$ to $V$, $g$ is a morphism
from $ V$ to $I$,we denote $f$ and  $g$ by
\[
\begin{tangle}
\Q f\\
\object{V}\\
\end{tangle}
\step \hbox { , } \step
\begin{tangle}
\object{V}\\
\QQ g\\
\end{tangle} \ \ \ \ .
\]

 If
$f$ is a morphism from $U \otimes V$ to $P$, $g$ is a morphism from
$U \otimes V$ to $I$ and $\zeta$  is a morphism from $U\otimes V$
to $V\otimes U$, we denote $f$, $g$ and $\zeta$ by
\[
\begin{tangle}
\object{U}\step[2]\object{V}\\
\tu f\\
\step\object{P}\\
\end{tangle}
\step,\step
\begin{tangle}
\object{U}\step[2]\object{V}\\
\coro g\\
\end{tangle}
\step,\step
\begin{tangle}
\object{U}\step[2]\object{V}\\
\ox \zeta\\
\object{V}\step[2]\object{U}\\
\end{tangle}
\step .
\]

   In particular, we denote the braiding
$C_{U,V}$  and its inverse $C_{U,V}^{-1}$ by
\[
\begin{tangle}
\object{U}\step[2]\object{V}\\
\x\\
\object{V}\step[2]\object{U}\\
\end{tangle}
\step,\step
\begin{tangle}
\object{V}\step[2]\object{U}\\
\xx\\
\object{U}\step[2]\object{V}\\
\end{tangle}
\step . \\
\]

$\xi $ is called  an $R$-matrix of ${\cal C}$ if
 $\xi$ is a natural
isomorphism from $\otimes $ to $\otimes \tau $ and  for any $U,V,W
\in ob {\cal C}$, the  Yang-Baxter equation of ${\cal C}$:

(YBE):
\[
\begin{tangle}
\object{U}\step[2]\object{V}\step[2]\object{W}\\
\ox \xi \step [2] \id\\
\id \step [2] \ox \xi \\
\ox \xi \step [2] \id \\
\object{W}\step [2]\object{V}\step[2]\object{U}\\
\end{tangle}
\ \ = \
\begin{tangle}
\object{U}\step[2]\object{V}\step[2]\object{W}\\
\id \step [2] \ox \xi \\
\ox \xi \step [2] \id\\
\id \step [2] \ox \xi \\
\object{W}\step [2]\object{V}\step[2]\object{U}\\
\end{tangle}\
\] holds.  In particular, the above equation is called the Yang-Baxter
equation on $V$ when $U = V=W.$

\begin {Lemma}\label {0.3}

(i) The braiding $C$ of braided tensor category  $ ({ \cal C}, C)$
is an $R$-matrix of $ { \cal C}$, i.e.
\[
\begin{tangle}
\object{U}\step[2]\object{V}\step[2]\object{W}\\
\x \step [2] \id\\
\id \step [2] \x \\
\x \step [2] \id \\
\object{W}\step [2]\object{V}\step[2]\object{U}\\
\end{tangle}
\ \ = \
\begin{tangle}
\object{U}\step[2]\object{V}\step[2]\object{W}\\
\id \step [2] \x \\
\x \step [2] \id\\
\id \step [2] \x \\
\object{W}\step [2]\object{V}\step[2]\object{U}\\
\end{tangle}\
\] holds.

(ii)
\[
\begin{tangle}
\object{U}\step[2]\object{V}\step[2]\object{W}\\
\id\step[2]\x\\
\x\step[2]\id\\
\id\step[2]\tu f\\
\object{W}\step[3]\object{P}\\
\end{tangle}
\step=\step
\begin{tangle}
\object{U}\step[2]\object{V}\step\object{W}\\
\tu f\step\id\\
\step\x\\
\step\object{W}\step[2]\object{P}\\
\end{tangle}
\step , \step
\begin{tangle}
\object{U}\step[2]\object{W}\step[2]\object{V}\\
\x\step[2]\id\\
\id\step[2]\tu f\\
\object{W}\step[3]\object{P}\\
\end{tangle}
\step=\step
\begin{tangle}
\object{U}\step[2]\object{W}\step\object{V}\\
\id\step[2]\hxx\\
\tu f\step\id\\
\step\x\\
\step\object{W}\step[2]\object{P}\\
\end{tangle} \ \ .
\]
\end {Lemma}

{\bf Proof.} (i)

\[
\begin{tangle}
\object{U}\step[2]\object{V}\step[2]\object{W}\\
\x \step [2] \id\\
\id \step [2] \x \\
\x \step [2] \id \\
\object{W}\step [2]\object{V}\step[2]\object{U}\\
\end{tangle} \ \ = \ \
\begin{tangle}

\step\object{U}\step[2]\object{V}\step[2]\object{W}\\
\step \id \step [2] \id \step [2] \id \\

\obox 3{C^{}_{V,W}}    \step [2] \id\\

 \step\id \step [2] \id \step [2] \id \\

\step\id \step [2] \obox 3 {C_{U,W}} \\

\step \id \step [2] \id \step [2] \id \\

\obox 3 {C_{V,W}} \step [2] \id \\

\step \id \step [2] \id \step [2] \id \\
\step\object{W}\step [2]\object{V}\step[2]\object{U}\\
\end{tangle} \ \ \stackrel {\mbox {by Hexagon Axiom }} {=}\begin{tangle}

\step\object{U \otimes V} \step[4]\object{W}\\
\step \id  \step [4] \id \\

\obox 3{C^{}_{V,W}}    \step [2] \id\\

 \step \id \step [3] \ne2 \step [2]\\

 \obox 4 {C_{V\otimes U,W}} \\

\step \id \step [2] \id \step [2] \\
\step\object{W}\step [3]\object{V \otimes U}\step[2]\\
\end{tangle} \ \
\]

\[\ \ \stackrel {\mbox {by naturality  }} {=}\begin{tangle}

\step\object{U \otimes V} \step[3]\object{W}\\
\step \id  \step [2] \ne1 \\

\obox 4 {C_{U\otimes V,W}}
   \step [2] \\

 \step \id \step [2] \id \step [2]\\

\step \id \step [1]\obox 3{C^{}_{U,V}}  \\

\step \id \step [2] \id \step [2] \\
\step\object{W}\step [3]\object{V \otimes U}\step[2]\\
\end{tangle} \ \
\stackrel {\mbox {by Hexagon Axiom }} {=}
\begin{tangle}
\object{U}\step[2]\object{V}\step[2]\object{W}\\
\id \step [2] \x \\
\x \step [2] \id\\
\id \step [2] \x \\
\object{W}\step [2]\object{V}\step[2]\object{U}\\
\end{tangle}\ \ \ .
\]

(ii)
 We have
\[
\begin{tangle}
\object{U}\step[2]\object{V}\step\object{W}\\
\tu f\step\id\\
\step\x\\
\step\object{W}\step[2]\object{P}\\
\end{tangle}
\step=\step
\begin{tangle}
\object{U\otimes V}\step[2]\object{W}\\
\O f\step[2]\id\\
\x\\
\object{W}\step[2]\object{P}\\
\end{tangle}
 \ \ \stackrel {\mbox {by naturality  }} {=} \step
\begin{tangle}
\object{U\otimes V}\step[2]\object{W}\\
\x\\
\id\step[2]\O f\\
\object{W}\step[2]\object{P}\\
\end{tangle} \]
\[
\stackrel {\mbox {by Hexagon Axiom }} {=}
\begin{tangle}
\object{U}\step[2]\object{V}\step[2]\object{W}\\
\id\step[2]\x\\
\x\step[2]\id\\
\id\step[2]\tu f\\
\object{W}\step[3]\object{P}\\
\end{tangle}\ \ \ .
\]
 Similarly, we can show the second equation. $\Box$

Now we give some concepts as follows: Assume  that $H, A \in ob \
{\cal C},$  and

\begin {eqnarray*}
\alpha : H \otimes A \rightarrow &A& , \hbox { \ \ \ \ }
\beta : H \otimes A \rightarrow H,    \\
\phi :  A \rightarrow H \otimes  &A&   , \hbox { \ \ \ \ }
\psi : H  \rightarrow H \otimes A,  \\
m _H : H \otimes H \rightarrow &H& , \hbox { \ \ \ \ }
m_A : A \otimes A \rightarrow A,    \\
\Delta _H :  H \rightarrow H \otimes  &H&   , \hbox { \ \ \ \ }
\Delta _A : A  \rightarrow A \otimes A, \\
\eta _H : I  \rightarrow &H& , \hbox { \ \ \ \ }
\eta _A : I \rightarrow A,    \\
\epsilon  _H :  H  \rightarrow   &I&   , \hbox { \ \ \ \ }
\epsilon  _A  : A  \rightarrow I.     \\
\end {eqnarray*}
are morphisms in ${\cal C}$.

$(A, m_A, \eta _A )$  is called an algebra living in ${\cal C}$, if
 the following conditions are
satisfied:
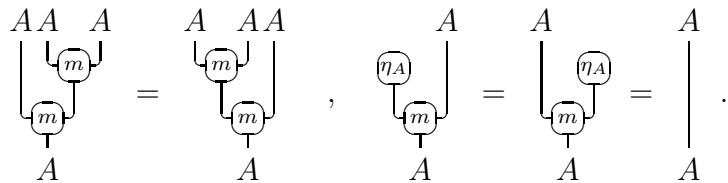
\[
\begin{tangle}
\object{A}\step\object{A}\step[2]\object{A}\\
\id\step\tu m\\
\tu m\\
\step\object{A}\\
\end{tangle}
\step=\step
\begin{tangle}
\object{A}\step[2]\object{A}\step\object{A}\\
\tu m\step\id\\
\step\tu m\\
\step[2]\object{A}\\
\end{tangle}
\step\step,\step\step
\begin{tangle}
\step[2]\object{A}\\
\Q {\eta_A}\step[2]\id\\
\tu m\\
\step\object{A}\\
\end{tangle}
\step=\step
\begin{tangle}
\object{A}\\
\id\step[2]\Q {\eta_A}\\
\tu m\\
\step\object{A}\\
\end{tangle}
\step=\step
\begin{tangle}
\object{A}\\
\id\\
\id\\
\object{A}\\
\end{tangle}
\quad .
\]
\noindent In this case, $\eta _A$ and $m_A$  are called  unit and
multiplication of $A$ respectively.

$(H, \Delta _H, \epsilon _H)$  is called a coalgebra living in
${\cal C},$ if the following conditions are satisfied:
\[
\begin{tangle}
\step\object{H}\\
\td \Delta\\
\id\step\td \Delta\\
\object{H}\step\object{H}\step[2]\object{H}\\
\end{tangle}
\step=\step
\begin{tangle}
\step[2]\object{H}\\
\step\td \Delta\\
\td \Delta\step\id\\
\object{H}\step[2]\object{H}\step\object{H}\\
\end{tangle}
\step\step,\step\step
\begin{tangle}
\step\object{H}\\
\td \Delta\\
\QQ {\varepsilon_H}\step[2]\id\\
\step[2]\object{H}\\
\end{tangle}
\step=\step
\begin{tangle}
\step\object{H}\\
\td \Delta\\
\id\step[2]\QQ {\varepsilon_H}\\
\object{H}\\
\end{tangle}
\step=\step
\begin{tangle}
\object{H}\\
\id\\
\id\\
\object{H}\\
\end{tangle}
\quad .
\]
 In this case, $\epsilon _H$ and $\Delta _H$  are called counit
and comultiplication of $H$ respectively.

If $A$  is an algebra and $H$  is a coalgebra, then $Hom_{\cal C}
(H, A)$  becomes  an algebra under the convolution product \[
f*g=\quad
\begin{tangle}
\step\object{H}\\
\td \Delta\\
\O f\step[2]\O g\\
\tu m\\
\step\object{A}\\
\end{tangle}\ \ \ .
\]
 \noindent and its unit element $\eta
= \eta _A \epsilon _H.$ If $S$ is  the inverse of $id _H$ in $Hom
_{\cal C} (H,H)$, then $S$ is called  antipode of $H$.

If $(H, m_H, \eta _H)$ is  an algebra,  and $(H, \Delta _H, \epsilon
_H)$ is a coalgebra living in ${\cal C}$, and the following
condition is satisfied:
 \[
\begin{tangle}
\object{H}\step[2]\object{H}\\
\tu m\\
\td \Delta\\
\object{H}\step[2]\object{H}\\
\end{tangle}
\step=\step
\begin{tangle}
\step\object{H}\step[3]\object{H}\\
\td \Delta\step\td \Delta\\
\id\step[2]\hx\step[2]\id\\
\tu m\step\tu m\\
\step\object{H}\step[3]\object{H}\\
\end{tangle}
\step,\step
\begin{tangle}
\object{H}\step[2]\object{H}\\
\tu m\\
\step\QQ \varepsilon\\
\end{tangle}
\step=\step
\begin{tangle}
\object{H}\step[2]\object{H}\\
\id\step[2]\id\\
\QQ \varepsilon\step[2]\QQ \varepsilon\\
\end{tangle}
\step,\step
\begin{tangle}
\step\Q \eta\\
\td \Delta\\
\end{tangle}
\step=\step
\begin{tangle}
\Q \eta\step[2]\Q \eta\\
\end{tangle}
\step,\step
\begin{tangle}
\Q \eta\\
\QQ \varepsilon\\
\end{tangle}
\step=\step
\begin{tangle}\object{I}\\
\id\\
\id\\
\object{I}
\end{tangle}
\step.\step
\]

\noindent then $H$ is called a bialgebra living in ${\cal C}$. If
$H$ is a bialgebra and there is  an inverse $S$ of $id _H$  under
convolution product in $Hom _{\cal C} (H, H)$, then $H$  is called a
Hopf algebra living in ${\cal C},$ or a braided Hopf algebra.

When $H$ is a Hopf algebra algebra in braided tensor category
$({\mathcal C}, C)$, then the condition above is equivalent to

(YD): \[ \phi \alpha =  \begin{tangle}
\step [3]\object{H}\step[7]\object{M} \\

\step \Cd \step [4] \td \phi \\

\cd \step [3] \O S \step [3] \ne2 \step [2] \id \\

\id \step [2] \nw1  \step [2]  \x \step [4] \id  \\

\id \step [2]  \step \x  \step[2] \id \step [4] \id\\

\id \step [2] \ne1  \step [2] \x \step[3]   \ne2\\

\cu  \step[2] \ne2  \step[2] \tu \alpha  \\

\step \cu  \step[4]  \step \id  \\
\step [2]\object{H}\step[6]\object{M}
  \end{tangle} \ \ \ \   .
  \]
  Let $^H_H {\mathcal YD (C)}$ denote the category of all  Yetter-Drinfeld
  $H$-modules in ${\mathcal C}$. If $({\mathcal C}, C) = {\mathcal V}ect(k)$, we write
  $^H_H {\mathcal YD (C)} = $ $^H_H{\mathcal YD}$, called Yetter-Drinfeld
  category.
  It follows from  \cite {RT93} and \cite [Theorem 4.1.1]{BD98} that
  $^H_H {\mathcal YD}({\mathcal C})$ is a braided tensor
  category with $^{ YD}C _{U,V} = (\alpha _V \otimes id_U) (id _H \otimes C_{U, V})
  (\phi _U\otimes id _V)$ for any
  two Yetter-Drinfeld modules  $(U, \phi _U, \alpha _U)$ and $(V, \phi _V, \alpha
  _V)$  when $H$ has an invertible antipode. In this case,
$^{YD}C _{U,V}^{-1} = (id _V \otimes \alpha _U)( C_{H,V}^{-1}
\otimes id _U) (S^{-1} \otimes C_{U, V}^{-1}) (C_{U,H}^{-1} \otimes
id _V)
  ( id _U \otimes \phi _V).$ Algebras, coalgebras and Hopf algebras
  and so on
  in $^H_H {\mathcal YD}$ are called  Yetter-Drinfeld ones or YD ones in
  short.

  Let $\mathbb Z,$ $\mathbb N$ and  $\mathbb C$ denote integer set, natural number set and complex field, respectively. Throughout basic field  is a complex  field $\mathbb C$, which is denoted by $k$ sometimes.
If $V$ is a vector space with a basis $v_1,  v_2,  \cdots,  v_n$  and  $q_{ij} \in \mathbb C^*$
for $1\le i,  j \le n$ such that map $C: \left \{
\begin
{array} {lll} V \otimes V& \rightarrow   & V\otimes V \\
 v_i \otimes v_j & \mapsto   & q_{ij} v_j \otimes v_i
\end {array} \right., $ then $(V,  C)$ is called a braided vector space of diagonal type.
Denote by $(q_{ij})_{n\times n}$ the braiding matrix of $(V,  C)$ under the
basis $v_1,  v_2,  \cdots,  v_n$. Then $(V, C)$ is also written as $(V, (q_{ij})_{n \times n}).$ We can get Nichols algebra $\mathfrak B(V)$  (see \cite {He06b, He06a} and \cite {ZZC04}).
 Let  $\mathcal B$ denote the set of all  generators  of restricted PBW basis and  $\mathcal P $ denote restricted PBW basis of $\mathfrak B(V)$.

Let  $v_1,  v_2,  \cdots,  v_n$ be a basis of $V\in ^{kG}_{kG} {\mathcal YD}$ with comodule operation and module operation  $\delta (v_i) = g_i \otimes v_i, $ $h \cdot v_i = \chi _i (h) v_i$, where $\chi _i$ is a multiplication character of $G$, i.e. a homomorphism from $G$ to $\mathbb C^*$, for $1\le i \le n$. It is clear that  $(V,  C)$ is a braided vector space of diagonal type. Otherwise, for any $q_{ij} \in \mathbb C$, $1\le i, j \le n$ and a basis $v_1,  v_2,  \cdots,  v_n$ of $V$, $V$ can become a Yetter-Drinfeld  module over group $\mathbb Z^n$ by defining $\delta (v_i ) = e_i \otimes v_i$ and $e_j \cdot v_i = q_ {ji}v_i$, $1\le i, j \le n$, where $e_1 = (1, 0, \cdots, 0), \cdots, e_n = (0, 0, \cdots, 1) \in \mathbb Z^n.$

 \section {Quantum Design of Industry}

In this section we give some examples to  design by means of braiding. That is, we  design by means of quantum theory.

\begin {Example} \label {2.2} (Category of quantum design grammars) Let $\{D_i \mid 1\le i \le m\}$ be the set of all shapes in some design $D$.
Let  $v_1,  v_2,  \cdots,  v_n$ be a basis of $V\in ^{kG}_{kG} {\mathcal YD}$ with $\delta (v_i) = g_i \otimes v_i, $ $h \cdot v_i = \chi _i (h) v_i$.
 Assume that   $g_i = g_j$ and  $\chi _i=  \chi _j$  when there exists rule  $R_{ij} : D_i \rightarrow D_j$. Define a morphism $f_{ij}$ from $kv_i$ to $kv_j$ in  $^{kG}_{kG} {\mathcal YD}$ such that $f_{ij} (v_i)=  R_{ij}(D_i)$. It is clear that  $f_{ij}$ is a morphism in   $^{kG}_{kG} {\mathcal YD}$ and $kv_i$ and  $V$, are in  $^{kG}_{kG} {\mathcal YD}$.  Define a map $\psi : \mathcal P \rightarrow  \{D_i \mid 1\le i \le m\}$ such that $\psi (v_i) = D_i$.  $\psi (u)$ denote the combination of shapes of  $D_{i_1}$, $D_{i_1}$, $\cdots,$ $D_{i_r}$ when $u = v_{i_1} v_{i_2} \cdots v_{i_r} \in \mathcal P.$

\end {Example}

\begin {Example} \label {2.3} ( Quantum design grammars of Coca-Cola)
Let $D_1, D_2$, $D_3$ and $ D_4$ be initial shapes; $D_5, D_6$, $D_7$ and $ D_8$ be cap shape, above shape, middle  shape and below shape of battle of Coca-Cola; Rule $R_{15} : D_1 \rightarrow D_5$, Rule $R_{26} : D_2 \rightarrow D_6$,  Rule $R_{37} : D_3 \rightarrow D_7$ and $R_{48} : D_4 \rightarrow D_8$.

(i) Let $G = \mathbb Z_7$ be cycle group  and $C(v_s\otimes v_t) = e^{\frac {2\pi i g_sg_t} {7}}= \omega ^{g_sg_t}$, where  $\omega := e^{\frac {2\pi i } {7}}  $ and $i := \sqrt{-1}$. Assume that     $g_s =s$, $\chi _s (g_t) = \omega ^{st}$ when $1\le s, t \le 4;$  $g_{4+p}= g_p$, $\chi _{4+p} = \chi _p$ for $1\le p \le 4.$

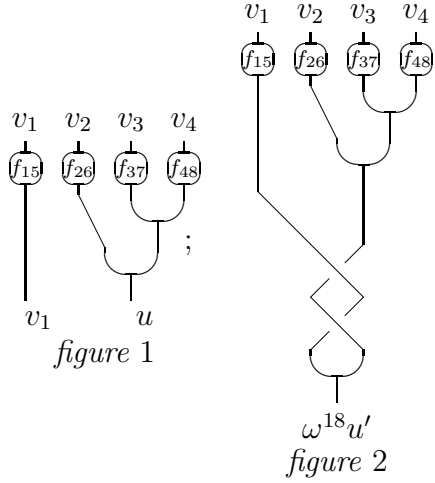
\[
\begin{tangles}{clr}
\step[1]\object{v_1}\step[2]\object{v_2}\step[2]\object{v_3}\step[2]\object{v_4}\\
\step[1]\O {f_{15}}\step[2]\O {f_{26}}\step[2]\O {f_{37}}\step[2]\O {f_{48}}\\
\step[1]\id \step[2]\nw1 \step[1]\cu \\
\id \step[3]\cu \\
\object{v_1}\step[4]\object{u}\\
\step[1]\object{\hbox {figure } 1}\\
\end{tangles} ;\ \ \
\begin{tangles}{clr}
\step[1]\object{v_1}\step[2]\object{v_2}\step[2]\object{v_3}\step[2]\object{v_4}\\
\step[1]\O {f_{15}}\step[2]\O {f_{26}}\step[2]\O {f_{37}}\step[2]\O {f_{48}}\\
\step[1]\id \step[2]\nw1 \step[1]\cu \\
\id \step[3]\cu \\
\nw2\step[3]\id \step[1]\\
\step[1]\x\\
\step[1]\x\\
\step [1]\cu\\
\step[1]\object{ \omega ^{18}u'}\\
\step[1]\object{\hbox {figure } 2} \\
\end{tangles} \ \ \ \ \  .
\]   Here $u= v_6v_7v_8$ and   $\psi (u)$ is the battle of Coca-Cola without cap, $u'= v_5v_6v_7v_8$ and   $\psi (u')$ is the battle of Coca-Cola with cap. Figure 1 is an ordinary design and Figure 2 is a quantum  design. Define that $\omega ^1 u', \omega ^2 u', \cdots, \omega ^7 u' $ means that the color of the battle of  Coca-Cola  is Red, orange, yellow, green, cyan, blue, purple, respectively. Therefore, $\omega ^{18}u' =\omega ^{4}u'$  and the color of the battle of  Coca-Cola  is green.

(ii) Let $G = \mathbb Z$ and $\chi _1(g_1) = q$; $1\not= q $ be a positive real number; $g_i =g_1$ and $\chi _i = \chi _1$  for $1\le i \le 8$. Then  $C(v_s\otimes v_t) = q(v_t \otimes v_s)$.
Similarly, we can get figure 3 by replacing $\omega ^{18}$ by $q^6$ in figure 2.
 Define that $q^6 u' $ means that the size  of the battle of  Coca-Cola  is  $q^6 $ times of  original  battle $ \psi (u')$ of  Coca-Cola.

\end {Example}

\begin {Example} \label {2.4} Let $D_0, D_1, D_3$, $D_5$,  $ D_7$ and $D_9$ be frame,  left steering wheel, left front gate, left back gate,  left front lamp and  left back lamp of car. We shall obtain a design of
 car by means of these left components of car as follows.
Let $G = \mathbb Z^{10}$ and  $C(v_0\otimes v_t) = qv_t\otimes v_0$, $C(v_s\otimes v_t) = v_t\otimes v_s$, where  $q\not=1$ is positive real number  for $s\not=0$, $s,  t = 1, 3, 5, 7, 9.$

\[
\begin{tangles}{lll}
\step[1]\object{v_9}\step[2]\object{v_7}\step[2]\object{v_5}\step[2]\object{v_3}   \step[2]\object{v_0}\step[2]\object{v_1}\step[2]\object{v_3}\step[2]\object{v_5}  \step[2]\object{v_7}   \step[2]\object{v_9}     \\

\step[1]\id \step[2]\id\step[2]\id \step[2]\x \step[2]\id\step[2]\id  \step[2]\id   \step[2]\id   \step[2]\id    \\

\step[1]\id \step[2]\id\step[2]\x  \step[2]\x\step[2]\id\step[2]\id\step[2]\id\step[2]\id\\

\step[1]\id \step[2]\x \step[2]\x \step[2]\cu\step[2]\id\step[2]\id\step[2]\id\\

\step[1]\x \step[2]\x\step[2] \id \step[2]\ne1\step[2]\ne1\step[1]\ne1\step[1]\ne1\\

\step[1]\id    \step[2]\x     \step[2]\id \step[2]\x\step[2]\ne1\step[1]\ne1\step[1]\ne1\\\\

\step[1]\cu\step[2]\id\step[2]\x\step[2]\cu  \step[1]\ne1\step[1]\ne1 \\

\step[2]\nw1\step[2]\x\step[2]\id \step[2]\ne1 \step[1]\ne1\step[1]\ne1 \\

\step[3]\cu\step[2]\id\step[2]\x\step[2]\id\step[2]\id\\

\step[4]\nw1  \step[2]\x\step[2]\cu\step[2]\id\\

\step[5]\cu\step[2]\id\step[2]\ne1\step[2]\ne1\\

\step[6]\nw1\step[2]\x\step[2]\ne1\\

\step[7]\cu\step[1]\ne1\step[1]\ne1\\

\step[8]\cu\step[1]\ne1\\

\step[9]\cu\\

\step[10]\object{ q^4u} \\

\step[10]\object{\hbox {figure } 4} \\

\end{tangles},
\]where $u =: v _0 v_1v_3^2v_5^2v_7^2v_9^2$ and $q^4u = v_0 v _1 (qv_3)v_3(qv_5)v_5 (qv_7)v_7 (qv_9)v_9$. Define that $(qv_1),$ $ (qv_3),$ $ (qv_5),$  $  (qv_7),$   $  (qv_9)$ denote right steering wheel, right front gate, right back gate,  right front lamp and  right back lamp of car, respectively.
Figure 4 is a braiding diagram which represents a combination with  left  components of car and whole car  consists of them. Notice that the steering wheel is in left side. If we require that  the steering wheel is in right side of car, we must
place the $v_1$ in the left hand side of $v_0$ in braiding diagram and have braiding $C(v_1 \otimes v_0).$

\end {Example}

\begin {Example} \label {2.5} Let $ D_4, D_1$, $D_2$ and  $D_3$  be   the first floor, the second floor and the third floor   of ship; $D_5$ the  brightness of ship. We shall obtain a design about the brightness of
 ship by means of braiding diagram as follows.
Let $G = \mathbb Z^5$ and  $C(v_s\otimes v_t) = q_{st}(v_t\otimes v_s)$, where  $q_{st}$ is positive real number for $1\le s, t  \le 5.$
\[
\begin{tangles}{lll}
\step[1]\object{v_3}\step[2]\object{v_2}\step[2]\object{v_1}\step[2]\object{v_4}  \step[2]\object{v_5}  \\

\step[1]\id \step[2]\id\step[2]\id \step[2]\x  \\

\step[1]\id \step[2]\id\step[2]\x  \step[2]\id  \\

\step[1]\id \step[2]\x \step[2] \id \step[2]\id\\

\step[1]\x \step[2]\id \step[2]\id\step[2]\id\\

\step[1]\object{v_5}\step[2]\object{u_3}\step[2]\object{u_2}\step[2]\object{u_1}  \step[2]\object{u_4}  \\

\step[5]\object{\hbox {figure 5}}  \\

\end{tangles}\ \ \ \ ,
\]where $ u_3:= q_{3 5}v_3$, $u_2:= q_{2 5}v_2$, $u_1:=q_{1 5}v_1$, $u_4:= q_{4 5}v_4.$
  Define that the brightness in   the first floor, the second floor,  the third floor and negative first  floor are  $q_{15}$ unit, $q_{25}$ unit, $q_{35}$ unit and $q_{45}$ unit according to $q_{3 5}v_1$,  $q_{2 5}v_2$,     $q_{3 5}v_2$ and $q_{4 5}v_4$ in braiding diagram. We can choose the value of $q_{st}$ according to the physical truth and requirement of customer.

\end {Example}
In fact, we also provide the brightness of every components of car  be means of braiding diagram.

\section {Quantum Design of Codes}
In this section we obtain a code by means of braiding.
\begin {Example} \label {3.1}We need send a word $w$ from $X$ to $Y$ by  internet and require to keep  secret, where $X$ and $Y$ are two persons or two companies. Let $\{D_i \mid 1\le i \le n \}$  be a vocabulary and  $Q = (q_{ij})_{n \times n}$   its quantum matrix.
\[
\begin{tangles}{lll}
\step[4]\object{v_1}\step[2]\object{v_2}\step[2]\object{v_3}\step[2]\object{v_4}   \step[2]\object{v_5}\step[2]\object{v_6}\step[2]\object{v_7}\step[2]\object{v_8}  \step[2]\object{v_9}   \step[2]\object{v_{10}}     \\

\object{  1}\step[4]\id \step[2]\id\step[2]\id \step[2]\x \step[2]\id\step[2]\id  \step[2]\id   \step[2]\id   \step[2]\id    \\

\object{2}\step[4]\id \step[2]\id\step[2]\x  \step[2]\x\step[2]\id\step[2]\id\step[2]\id\step[2]\id\\

\object{3}\step[4]\id \step[2]\x \step[2]\x \step[2]\cu\step[2]\id\step[2]\id\step[2]\id\\

\object{4}\step[4]\x \step[2]\x\step[2] \id \step[2]\ne1\step[2]\ne1\step[1]\ne1\step[1]\ne1\\

\object{5}\step[4]\id    \step[2]\x     \step[2]\id \step[2]\x\step[2]\ne1\step[1]\ne1\step[1]\ne1\\\\

\object{6}\step[4]\cu\step[2]\id\step[2]\x\step[2]\cu  \step[1]\ne1\step[1]\ne1 \\

\object{7}\step[5]\nw1\step[2]\x\step[2]\id \step[2]\ne1 \step[1]\ne1\step[1]\ne1 \\

\object{8}\step[6]\cu\step[2]\id\step[2]\x\step[2]\id\step[2]\id\\

\object{9}\step[7]\nw1  \step[2]\x\step[2]\cu\step[2]\id\\

\object{10}\step[8]\cu\step[2]\id\step[2]\ne1\step[2]\ne1\\

\object{11}\step[9]\nw1\step[2]\x\step[2]\ne1\\

\object{12}\step[10]\cu\step[1]\ne1\step[1]\ne1\\

\object{13}\step[11]\cu\step[1]\ne1\\

\object{14}\step[12]\cu\\

\step[13]\object{ qu} \\

\step[12]\object{\hbox {figure 6}}  \\

\end{tangles}\ \ \ \ ,
\] where $u= v_5v_6v_{4} v_7 v_3 v_8 v_2 v_9 v_1 v_{10}$.

Let
$Q$ and figure 6 are the private key cryptography. i.e. $X$ and $Y$ have $Q$ and figure 6 without using internet. Assume  $w = D_s$. We choose $p$ such that  $s \in \{p, p+1, \cdots, p+9\}$ and positive real numbers $a_1, \cdots, a_{10}$ randomly with $a_1 \not=1, \cdots, a_{10}\not=1$. Let $b_i = a_i$ when $i\not= s$ and $b_s =1.$   Let $v_i := D_{p+i-1}$ and $\widetilde{ q_{ij}} := q_{p+ i-1, p+j-1}$ for $1\le i, j \le 10.$ It is clear $C^{-1} (v_i \otimes v_j) = \widetilde{q_{ji}}^{-1}(v_j \otimes v_i)$ and $C^{-1} (v_i \otimes u) =  \prod _{j =1} ^{10} \widetilde{q_{ji}}^{-1} (u\otimes v_i)  $.  Compute the value  $c_i$ of $b_i \prod _{j =1} ^{10} \widetilde{q_{ji}} ^{-1}$  and send $c_i$  to $Y$ for
$i =1, 2, \cdots, 10$ and send $p$ to $Y$. In $Y$  compute the value  $d_i$ of $c_i \prod _{j =1} ^{10} \widetilde{q_{ij}} $ and $d_i = 1$ if and only if $s=i+p -1$ since   coefficient  of  $C^{-1} (v_i \otimes u)$ is $\prod _{j =1} ^{10} \widetilde{q_{ji}}^{-1}$ and coefficient  of  $c (u\otimes v_i)$ is $\prod _{j =1} ^{10}\widetilde{ q_{ji}}$. Consequently, we can get the word $w$ in $Y.$  Furthermore, we can send a article from $X$ to $Y$   by means of internet and  keeping  secret because they consist of some single words. We also can send file  from $X$ to $Y$   by means of internet and  keeping  secret because the password  consist of some single words.

\end {Example}

Otherwise, we can modify the method above as follows.

Remark. (i)  In case above, we choose the last layer 14-th layer and left line $u$. In fact, we also can choose $t$-th layer with $1\le t\le 14$ and left line $v$ which need  be sent to $Y$. That is,
 private key cryptography contains $Q$,  figure 6, $p$,  $t$-th layer  and left line $v$.
For example, set $t= 7$, i.e. $v = v_5v_6.$
It is clear  $C^{-1} (v_i \otimes v) =  \prod _{j =5} ^{6} \widetilde{q_{ji}}^{-1} (u\otimes v_i)  $ and $c (v \otimes v _i) =  \prod _{j =5} ^{6} \widetilde{q_{ji} }(v _i\otimes v)$. Compute the value  $c_i$ of $b_i \prod _{j =5} ^{6} \widetilde{q_{ji}} ^{-1}$  and send $c_i$  to $Y$ for
$i =1, 2, \cdots, 10$. In $Y$  compute the value  $d_i$ of $c_i \prod _{j =5} ^{6} \widetilde{q_{ij}} $ and $d_i = 1$ if and only if $s=i$.

(ii)  figure 6 and  vocabulary $\{D_i \mid 1\le i \le n \}$  can become public key cryptography. The order of vocabulary $\{D_i \mid 1\le i \le n \}$ can be the same as  Xinhua dictionary.

(iii) Vocabulary $\{D_{r+i} \mid 1\le i \le n \}$  can become public key cryptography and $r$ can become the private  key cryptography.

\end {document}